\documentclass[preprint,12pt,authoryear]{elsarticle}
\usepackage{graphicx,lmodern,amsfonts,mathrsfs,amsmath,amssymb,latexsym,xspace,epsfig,psfrag,color,tikz,tikzscale,pgfplots,colortbl}
\usepackage[flushleft]{threeparttable}
\usepackage[colorlinks,bookmarksnumbered,citecolor=blue]{hyperref}
\usepackage[ruled]{algorithm2e}
\usepackage[utf8]{inputenc}  
\usepackage{scalefnt}
\usepackage{subfigure}
\usepackage{multirow}
\usepackage[flushleft]{threeparttable}
\usepackage{pdflscape}
\usepackage{multicol}
\usepackage{natbib}
\usepackage{hyperref}
\usepackage{lineno}
\usepackage{longtable, lscape} 
\usepackage{bbm}
\usepackage{calrsfs}
\usepackage{txfonts}
\usepackage{rotating}
\usepackage{comment}

\usepackage[top=30mm,bottom=20mm,left=30mm,right=20mm,pdftex,includeheadfoot]{geometry}

\usepackage{booktabs}
\usepackage{multirow}
\usepackage{pgfplotstable}
\usepackage{filecontents}
\usepackage{xspace}
\usepackage{tabularx}
\usepackage{siunitx}

\usepackage{tikz,pgfplots}
\pgfplotsset{compat=1.15}
\usepgfplotslibrary{statistics}
\usetikzlibrary{shapes,arrows}

\definecolor{Gray}{gray}{0.95}
\definecolor{Gray1}{gray}{0.6}
\definecolor{Gray2}{gray}{0.8}
\definecolor{Gray3}{gray}{1}

\journal{Computers \& Operations Research}

\begin{document}

\begin{frontmatter}

\title{An Adaptive Iterated Local Search Heuristic for the Heterogeneous Fleet Vehicle Routing Problem}

\author[unifesp]{Vinícius R. Máximo}
\ead{vinymax10@gmail.com}

\author[montreal]{Jean-François Cordeau}
\ead{jean-francois.cordeau@hec.ca}

\author[unifesp]{Mariá C. V. Nascimento\corref{cor1}}
\ead{mcv.nascimento@unifesp.br}
\cortext[cor1]{Corresponding author}

\address[unifesp]{Instituto de Ciência e Tecnologia, Universidade Federal de São Paulo (UNIFESP)\\ Av. Cesare M. G. Lattes, 1201, Eugênio de Mello, São José dos Campos-SP, CEP: 12247-014, Brazil
}

\address[montreal]{HEC Montréal and GERAD, 3000 chemin de la Côte-Sainte-Catherine, Montréal, H3T 2A7, Canada
}

\begin{abstract}
The Heterogeneous Fleet Vehicle Routing Problem (HFVRP) is an important variant of the classical Capacitated Vehicle Routing Problem (CVRP) that aims to find routes that minimize the total traveling cost of a heterogeneous fleet of vehicles. This problem is of great interest given its importance in many industrial and commercial applications. In this paper, we present an Adaptive Iterated Local Search (AILS) heuristic for the HFVRP. AILS is a local search-based meta-heuristic that achieved good results for the CVRP. The main characteristic of AILS is its adaptive behavior that allows the adjustment of the diversity control of the solutions explored during the search process. The proposed AILS for the HFVRP was tested on benchmark instances containing up to 360 customers. The results of computational experiments indicate that AILS outperformed state-of-the-art metaheuristics on 87\% of the instances.

\end{abstract}


\end{frontmatter}


\section{Introduction}


The Vehicle Routing Problem (VRP) is one of the most studied operations research problems in the literature \citep{Penna2013}. Its classical form was introduced by \cite{Dantzig1959} and consists of finding minimum traveling cost routes for a homogeneous fleet of vehicles while respecting the capacity of these vehicles. In most practical applications with a large number of customers, heuristics, in particular metaheuristics, are more appropriate than exact methods for solving VRPs in reasonable time. The hill-climbing and improvement mechanisms in metaheuristics allow them to achieve excellent performance in many of the most important problem variants.
 

The literature contains many variants of the VRP that address specific characteristics of practical applications. In particular, those involving vehicles of the same fleet with different features are of great interest in practice. This study focuses on this well-known VRP variant called the Heterogeneous Fleet Vehicle Routing Problem (HFVRP). The HFVRP is a hard-to-solve problem since, in addition to deciding the composition of the routes, it is necessary to choose the appropriate vehicle for each route. While some authors have introduced exact solution methods for the HFVRP \citep{Baldacci2010exact,Baldacci2008_2}, the vast majority of the HFVRP studies consider heuristic methods \citep{Liu2009,Li2007}.

Iterated Local Search (ILS) is a metaheuristic that was successfully applied to the HFVRP by \citet{Penna2013} and \citet{Subramanian2012}. ILS is composed of two stages: perturbation and local search. Initially, the construction of an initial solution provides a starting point, and the perturbation and local search operators are iteratively applied throughout the search process. After performing local search, an acceptance criterion decides whether to use the current solution in the following iterations as a reference solution. Despite the high-quality solutions ILS obtains on average, the results obtained on the HFVRP indicate that there are instances for which the method gets stuck in local optima. Investigating the reasons behind this behavior, empirical tests pointed out that the difficulty of ILS to escape from local optima is often due to the criterion employed for ILS to accept solutions (acceptance criterion) and the perturbation strength (degree of perturbation). The choice of the acceptance criterion and the degree of perturbation are here called the diversity control problem.  The related diversity parameters need to be carefully adjusted, especially in large routing instances, for ILS to achieve a good performance.  

\cite{Maximo2021} proposed an Adaptive Iterated Local Search (AILS) applied to the  Capacitated Vehicle Routing Problem (CVRP) that implements diversity controls both in the degree of perturbation and in the acceptance criterion. This method achieved greater robustness in comparison to other algorithms proposed in the literature. Inspired by this previous work, in this paper, we propose an AILS for the HFVRP. We carried out computational experiments using benchmark instances ranging from 20 to 360 vertices. The AILS was compared to other algorithms from the literature such as the ones proposed by \citet{Penna2013},  \citet{Subramanian2012} and  \citet{Vidal2014}. The results indicate that AILS outperformed the other solution methods in 87\% of the instances.



The rest of the paper is organized as follows. Section \ref{sec:problem} defines the HFVRP and the main notation used throughout the paper. Section \ref{sec:related} reviews the main existing metaheuristics for the HFVRP. Section \ref{sec:proposedAILS} describes the proposed AILS to solve the HFVRP. In Section \ref{section:EC}, the computational experiments are presented as well as the results of the proposed AILS in comparison to the state-of-the-art heuristics. Finally, in Section \ref{section:conclusao}, we conclude the paper with some final remarks and directions for future work.

\section{Problem description}\label{sec:problem}

The HFVRP, is defined as follows: Let $G=(V,E)$ be a complete undirected graph where $V = \{0,1,2,\ldots,n\}$ is the set of $n+1$ vertices and $E = \{ (i,j) : i,j \in V, i  <  j\}$ is the set of edges. The depot is represented by vertex 0 and the set $V_c = V\setminus{\{0\}} $ represents the set of customers that must be visited once. Each edge $(i,j) \in E$ has a length $d_{ij}$ which represents the distance from vertex $i$ to vertex $j$. Each customer $i \in V_c$ has a non-negative demand $q_i$ that must be satisfied. The vehicle fleet is represented by the set $T = \{1,2,\ldots,h\}$ and there are thus $h$ types of vehicles available.

Each route $i$ will be represented by the pair $(R_i,t_i)$  where $R_i = \{v^{i}_0, v^{i}_1,\ldots,v^{i}_{m_i}\}$ represents an ordered set of $m_i +1$  vertices, where $v^{i}_0=v^{i}_{m_i}=0$ and $t_i \in T$ represents the vehicle type assigned to route $i$. The route is a closed loop that contains the depot and no repetition of vertices. Each cycle is represented by a sequence of vertices so that one vertex will be adjacent to another vertex if they are consecutive in the sequence. The routes of a solution $s$ are represented by the set $\mathscr{R}^s=\{(R^s_1,t^s_1),(R^s_2,t^s_2),\ldots, (R^s_{m_s},t^s_{m_s})\}$ where $m_s$ indicates the number of routes in solution $s$, $R^s_i \cap R^s_j =\{0\}$, for $i\neq j$, and $\cup_{i=1}^{m_s} R^s_i = V$.

Each type of vehicle $t \in T$ has a capacity $\bar{q}_t$, a unit cost $r_t$, a fixed operating cost $f_t$ and number $u_t$ of available vehicles. The capacity $\bar{q}_t$ limits  the total demand that a vehicle can serve. Therefore, constraint $\sum_{v \in R_i} q_v \leq \bar{q}_t$ must be satisfied for every route $R_i$. The unit cost $r_t$ is a factor that multiplies the distance so that the travel cost of edge $(i,j)$ is given by $c_{ij}=d_{ij} r_t$. The fixed operating cost $f_t$ is a value that will be added to the objective function whenever a vehicle of type $t$ is used. The HFVRP belongs to the $\mathcal{NP}$-hard problem class as it is a generalization of the CVRP which has identical vehicles.

The existing variants of the HFVRP differ in the characteristics regarding the fleet size, which can be limited or not, and the costs involved in using a vehicle, which can be fixed or dependent on the distance traveled. This paper deals with five variants defined by \cite{Penna2013}: the HFVRP with fixed and dependent costs (HVRPFD); the HFVRP with limited fleet and only dependent costs (HVRPD); the FSM with unlimited fleet, fixed and dependent costs  (FSMFD); the FSM with unlimited fleet and fixed costs (FSMF); and the FSM with unlimited fleet and dependent costs (FSMD). Table \ref{variantes} highlights their main characteristics.

\begin{table}[!ht]
\caption{Characteristics of the five variants.}
\label{variantes}
\centering
\begin{tabular}{cccc}
\hline
\hline
Variant & Limited fleet & Fixed costs & Dependent costs \\
\rowcolor{Gray}
HVRPFD & \checkmark &\checkmark&\checkmark\\
HVRPD &\checkmark&&\checkmark\\
\rowcolor{Gray}
FSMFD &&\checkmark&\checkmark\\
FSMF &&\checkmark&\\
\rowcolor{Gray}
FSMD &&&\checkmark\\
\hline
\hline
\end{tabular}
\end{table}


\section{Related Works}\label{sec:related}



The HFVRP is a generalization of the classical capacited vehicle routing problem, the so-called CVRP. The main difference between the HFVRP and the CVRP is that the latter considers a fleet of identical vehicles. As a consequence, a solution method designed for the HFVRP can also solve the CVRP. The current state-of-the-art heuristics for the CVRP are POPMUSIC \citep{Queiroga2021}, 
Slack Induction by String Removals (SISRs) \citep{Christiaens2020} and AILS-PR \citep{Maximo2021}, which is AILS hybridized with path-relinking.

The HFVRP is a problem largely studied in the literature due to its many practical applications. In real routing problems, heterogeneous fleets are more frequent than homogeneous ones \citep{Hoff2010}. \cite{Baldacci2008_3} presented a survey of the HFVRP discussing the best algorithms proposed at the time and some variants of this problem. More recently, \cite{Koc2016} reviewed the last 30 years of literature on the HFVRP. The authors provided an overview of the evolution of meta-heuristics proposed for this problem. Since the literature is quite extensive, we limit the scope of this section to the state-of-the-art meta-heuristics for this problem.


\cite{Penna2013} proposed the first ILS for the HFVRP. The authors presented an ILS with Variable Neighborhood Descent and a random neighborhood ordering (RVND) in the local search. They named the method ILS-RVND. The local search uses the widely employed moves between routes known as 2-interchange \citep{Osman1993} and 2-opt* \citep{Potvin1995}. The perturbation process considers vertex swap movements, vertex shifts, and route splitting. The experiments were performed on  instances proposed by \cite{Golden1984} whose number of vertices varies from 20 to 100. The results of ILS-RVND outperformed previous heuristic methods introduced by \cite{Tarantilis2004}, \cite{Li2007}, \cite{Prins2009} and \cite{Brandao2011}.

\cite{Subramanian2012} proposed a hybrid algorithm composed by the ILS proposed by \cite{Penna2013} and a Set Partitioning (SP) formulation. The columms in the SP model correspond to routes generated by ILS and the formulation is solved by a Mixed Integer Programming (MIP) solver. The resulting algorithm, called ILS-RVND-SP, was tested on the instances introduced by \cite{Golden1984}, \cite{ Brandao2011} and \cite{Li2007} with up to 360 vertices. ILS-RVND-SP outperformed the ILS-RVND proposed by \cite{Penna2013} and found a new best solution to eight instances of the problem. The ILS-RVND-SP was significantly faster than the ILS-RVND. 

\cite{Vidal2014} presented a genetic algorithm called Unified Hybrid Genetic Search (UHGS) that solved numerous VRP variants, including three HFVRP variants (FSMFD, FMSF, FSMD). UHGS consists of a local search that explores a neighborhood defined by the 2-interchange and 2-opt moves. The solutions generated by the genetic algorithm are stored as a single sequence, that is, without the separation of routes. To generate a feasible solution, the authors use a splitting procedure proposed by \citet{Prins2004} that obtains the optimal separation of the routes. The authors allow infeasible solutions during the search through objective function penalties. This strategy, according to the authors, allowed a significant gain in the performance of the method. In the computational experiments considering the three HFVRP variants, they applied their solution method to the instances proposed by \cite{Golden1984}. The authors compared the results obtained by UHGS with those achieved by ILS-RVND-SP. UHGS found better solutions in most of the considered instances, despite the significantly higher computational times. UHGS presented better average results than ILS-RVND-SP on the vast majority of the tested instances. However, in comparison to ILS-RVND-SP, the average times of UHGS were approximately nine times longer.

\section{Adaptive Iterated Local Search (AILS)}\label{sec:proposedAILS}

Adaptive Iterated Local Search (AILS) \citep{Maximo2021} is a variant of ILS, a metaheuristic framework first introduced by \cite{Lourenco2003}. ILS is composed of two stages: perturbation and local search. These stages iteratively repeat until a stopping criterion is met. The primary difference between ILS and AILS is that the latter considers adaptive mechanisms to improve the search process. \citet{Maximo2021} proposed a diversity control that acts both on the perturbation degree and the acceptance criterion. 

The perturbation stage of AILS is composed of removal heuristics that define which vertices must be removed from routes of a given solution to be inserted in other positions of the solution according to insertion heuristics. The number of vertices to be removed is called the perturbation degree and is a key parameter of AILS. Controlling the perturbation degree allows the algorithm to have a constant and equivalent diversity for each vertex removal heuristic during the search process. This feature allows the algorithm to self-adjust, correcting the diversity degree during execution. The acceptance criterion is also a crucial mechanism to enable the algorithm to escape from local optima. In AILS for the CVRP, the authors empirically demonstrate that the more relaxed, the more likely it is to escape from local optima. Therefore, they adjust this criterion bearing this behavior in mind. In this paper, we adapt the AILS to the HFVRP. Algorithm \ref{AILS} describes the proposed method.

\begin{algorithm}[!htb]
\LinesNumbered
\SetAlgoLined
\KwData{Instance data}
\KwResult{The best solution found $s^*$}
$s \leftarrow$ Construct an initial solution \\ 
$s^r,s^* \leftarrow$ Local Search($s$) \\
\Repeat{stop criterion is met}
{ 
    \Repeat{ $s$ is feasible }
    {
        Choose at random a removal heuristic $ \mathscr{H}_k \in \mathscr{H}$\\
        $s \leftarrow$ Perturbation Procedure  ($s^r$, $\mathscr{H}_k$)\\
        $s \leftarrow$ Feasibility($s$)\label{linha5ails}\\
    }
   
    $s \leftarrow$ Local Search($s$)\label{linha5ails}\\
	Update the diversity control parameter $\omega_{\mathscr{H}_k}$ considering the distance between $s$ and $s^r$\\
    $s^r\leftarrow $ Apply acceptation criterion to $s$\\
	Update the acceptance criterion\\
	Assign $s$ to $s^*$ if $ f(s) < f(s^*)$\\
}
\caption{Adaptive Iterative Local Search}
\label{AILS}
\end{algorithm}

First, a procedure, described in Section~\ref{sec:solucaoInicial}, finds a feasible initial solution $s$ for the target problem. Then, the local search, introduced in Section~\ref{FacBL}, is applied to $s$ and the resulting local optimum is assigned as reference solution $s^r$ and best obtained solution $s^*$. Then, the main loop of the algorithm starts. In line 5 of  Algorithm \ref{AILS},  a vertex removal heuristic $\mathscr{H}_k$ to be applied to the reference solution $s^r$ in the perturbation procedure is randomly selected.   The perturbation strategy, described in Section \ref{perturbacao}, returns a solution $s$ that can violate the capacity constraints. Therefore, we apply the feasibility method, described in Section~\ref{FacBL}, to $s$ to attempt (with no guarantee) to obtain a  feasible solution in the neighborhood of $s$. If the resulting solution is feasible, a local search is applied to $s$; otherwise, the method returns to line 5. After the local search,  the diversity control adjustment step updates the diversity control parameter $\omega_{\mathscr{H}_k}$  that is linked to the selected removal heuristic $\mathscr{H}_k$. For this, we analyze the distance between $s$ and the reference solution $s^r$ to adjust the degree of perturbation $\omega_{\mathscr{H}_k}$.  This adjustment is described in Section \ref{ControlePert}. Then, if $s$ meets the acceptance criterion, described in Section \ref{criterioAceitacao}, it replaces the reference solution. Finally, the best overall solution is updated if the objective function value of $s$ is better than $s^*$, i.e. $f(s)<f(s^*)$.  Then, the stopping criterion -- presented in Section~\ref{parametertuning} -- is checked and if it has not been reached, the algorithm continues the search from the reference solution $s^r$.

\subsection{Construction of the Initial Solution}
\label{sec:solucaoInicial}

In this section, we describe the method for constructing an initial solution for the HFVRP. Algorithm \ref{alg:construcaoinicial} presents the pseudocode of the strategy used in this paper.
 
 \begin{algorithm}[!htb]
 \caption{Initial Solution}
\label{alg:construcaoinicial}
\LinesNumbered
\SetAlgoLined
\KwData{Data Instance}
\KwResult{A solution $s$}
\Repeat{ $s$ is feasible }
{
$V_a\leftarrow V_c$\\
$\mathscr{R}^s \leftarrow \emptyset$\\
\For{$j=1,\ldots, \underline{m}$}{
 Construct $R_j=\{v_0^{j},v_1^{j}, v_0^{j}\}$ by randomly choosing $v_1^j$ from $V_a$\\
 $ (R_j,t) \leftarrow $ Include in $R_j$ an available type of vehicle $t \in T$.\\
 $\mathscr{R}^s \leftarrow  \mathscr{R}^s  \cup (R_j,t)$\\
 $V_a \leftarrow V_a\backslash \{v_1\}$\\
} 
\While{$V_a\neq \emptyset$}
{ 
    Pick $v\in V_a$ at random\\
    Add  $v$ to solution $s$ following the insertion heuristic $ \mathscr{A}_1$\\
    $V_a \leftarrow V_a \backslash  \{v\} $\\
}
$s \leftarrow$ Feasibility($s$)\label{linha5ails}\\
}
\Return{$s$}
\end{algorithm}
 
In Algorithm~\ref{alg:construcaoinicial}, let $V_a$ be a set containing all  the vertices to be inserted into solution $s$, that is initially the whole set of customers. Then, the algorithm initializes the set of routes $\mathscr{R}^s$ as an empty set. Then, the process of constructing $\underline{m}$ routes with a single vertex occurs in lines 4 to 9. The value of $\underline{m}$ represents the minimum number of routes necessary to generate a feasible solution for the HFVRP, calculated according to Equation \eqref{minferior}:
 \begin{align}
&\underline{m}=\left\lfloor \frac{1}{\max_{t \in T} \bar{q}_t} \sum_{i \in V_c} q_i \right\rfloor.& \label{minferior}
\end{align}

For constructing a route $R_j$ besides randomly choosing a vertex from $V_a$, an available type of vehicle $t\in T$ is selected. It is worth mentioning that a type of vehicle $t$ is available if the number of vehicles of type $t$ already added to the solution is lower than $u_t$. 

From lines 10 to 14, the vertices not chosen to be added to the $\underline{m}$ routes are individually inserted to the solution by the insertion heuristic $\mathscr{A}_1$ described in Section \ref{perturbacao}. The resulting solution may not be feasible, given that the construction process did not take the vehicle capacity assigned to each route into account. Hence, a feasibility procedure is applied to solution $s$. However, the feasibility method does not guarantee the return of a feasible solution. For this reason, a new initial solution is built until the feasibility method can achieve a feasible solution.



\subsection{Perturbation}
\label{perturbacao}

The perturbation process produces changes in a solution $s$ following a heuristic to obtain a neighboring solution $s'$ potentially different from $s$. The perturbation methods adopted in this paper do not guarantee the feasibility of the resulting solution. For this reason, we use a feasibility algorithm described in Section \ref{FacBL} that will address this issue. The perturbation process has two types of heuristics: the removal and insertion heuristics. The removal heuristic decides which vertices to remove from the current solution.  The insertion heuristics determines where to place the removed vertices into the solution. In this paper, we adopt three removal heuristics and two insertion heuristics.

\textbf{Concentric removal, $\mathscr{H}_1$: } This strategy removes a randomly chosen vertex $v_r$ and the $\omega_{\mathscr{H}_1}-1$ closest vertices to $v_r$. It is called concentric because the removed vertices are within a covering radius where $v_r$ is the center.

\textbf{Random removal, $\mathscr{H}_2$: } This heuristic consists of removing  $\omega_{\mathscr{H}_2}$ randomly chosen vertices from the solution.

\textbf{Removal of a sequence of vertices, $\mathscr{H}_3$: } This heuristic removes from the solution $\omega_{\mathscr{H}_3}$ adjacent vertices belonging to a route. The starting vertex of the sequence is chosen at random. If $\omega_{\mathscr{H}_3}$ is greater than the size of the route, then the heuristic deletes this route and picks another route to remove the remaining vertices --   until  $\omega_{\mathscr{H}_3} $ vertices are removed.

The insertion heuristics  $\mathscr{A}_1$ and $\mathscr{A}_2$ do not allow a removed vertex $v_i$ to be adjacent to any of the vertices to which they were previously adjacent. Without loss of generality, let  $v_{i-1}$ and  $v_{i+1}$ be the vertices to which $v_i$ was adjacent. Then, the new adjacent vertices, let us say $v_j$ and $v_{j+1}$ necessarily meet the following conditions: $j\neq i-1$, $j\neq i+1$, $j+1\neq i$ and $j+1\neq i-1$.  Thus, we define the set $V_p(v_i)$ as the set of vertices  that respect these conditions, i.e., is the set of vertices that are candidates to be adjacent to $v_i$.


\textbf{Insertion by distance, $\mathscr{A}_1$:} This heuristic consists of inserting a vertex $v_i$ between vertices $v_{\hat{j}}$ and  $v_{\hat{j}+1}$,  such that $v_{\hat{j }} = \arg \min_{j \in V_p(v_i)} d_{i,j}$.


\textbf{Insertion by cost, $\mathscr{A}_2$:} This heuristic consists of inserting a vertex $v_i$ adjacent to vertices $v_{\hat{j}}$ and  $v_{\hat{j}+1}$, such that $v_{\hat{j} } = \arg \min_{j \in V_p(v_i)} \{c_{i,j}+c_{i,j+1} - c_{j,j+1}\}$.

Algorithm~\ref{pertubacao} presents the pseudocode of the introduced perturbation strategy.

\begin{algorithm}[!htb]
\LinesNumbered
\SetAlgoLined
\KwData{Solution $s^r$ and $\mathscr{H}_l$}
\KwResult{Solution $s'$}
$s' \leftarrow s^r$\\
Change by 1 unit the number of routes $m$ of the solution $s'$  with probability $\alpha$ if the new value of $m \in [\underline{m},\overline{m}]$\\
Randomly select two routes in solution $s'$ and with probability $\alpha$ change their type of vehicle \\
Randomly select an insertion heuristic $\mathscr{A}_j \in \mathscr{A}$\\
\For{$k=1,\ldots, \omega_{\mathscr{H}_l}$}
{
    Remove a vertex $v_i$ from $s'$ following heuristic \textbf{ $\mathscr{H}_l$}\\
    Insert $v_i$ in a route chosen according heuristic $\mathscr{A}_j$\\
}
\caption{Perturbation $\mathscr{H}_l$}
\label{pertubacao}
\end{algorithm}

Algorithm~\ref{pertubacao} initially makes a copy of the reference solution $s^r$  in $s'$. Then, the number of routes of the solution $s'$ is defined. It can increase or decrease by one unit with respect to the input solution $s$ with probability $\alpha$ as long as $m \in [\underline{m},\overline{m}]$, where

\begin{align}
\overline{m} = \left \{ \begin{array}{ll}
\ n, & \textrm{for FSMFD, FSMF and FSMD variants}  \\
 \sum_{t \in T} u_t, & \textrm{ for HVRPFD and HVRPD variants.} 
\end{array}\right. & \label{moverline}
\end{align}

The decision between increasing and decreasing the number of routes by one unit is random in the case the algorithm chooses to change the number of routes since it is not a bound of the interval (if $m$ is $\underline{m}$ or $\overline{m}$, the number of routes can only increase or decrease, respectively).  A randomly chosen route must be removed from $s'$ if the choice is to decrease the number of routes. In this case, the vertices of the removed route are inserted into other routes through the  insertion heuristic $\mathscr{A}_j$. If $m$ is increased by 1, then a new empty route will be included in the solution.

Then, the algorithm changes the type of vehicles of a pair of randomly selected routes, again, with probability $\alpha$. For HVRPFD and HVRPD variants, the type of vehicles of the two routes is swapped. For the FSMFD, FSMF and FSMD variants,  their vehicles are exchanged for other vehicles chosen according to a probability distribution proportional to the types of vehicles already present in solution $s$, as defined in Equation \eqref{prob}: 



\begin{align}
P(t,s) = \frac{|\{ (R_j,t') \in \mathscr{R}^s ~|~ t'=t \}| + 1}{m+h}.&\label{prob}
\end{align}

After this step, an insertion heuristic $\mathscr{A}_j \in \mathscr{A}$ is selected randomly. Then, the process of removing and adding $\omega_{\mathscr{H}_l}$ vertices of the solution $s'$ occurs for $\omega_{\mathscr{H}_l}$ iterations.
 
\subsection{Feasibility and Local Search}
\label{FacBL}

Algorithms \ref{alg:construcaoinicial} and \ref{pertubacao} produce solutions that may violate the capacity constraint of the vehicles assigned to each route. For this reason, a strategy that searches for feasible solutions in the neighborhood of infeasible ones might be necessary. Moreover, a local search method is applied to feasible solutions to return local optima. The feasibility and local search used in this paper are the same presented in \citet{Maximo2021}.

Both feasibility and local search neighborhoods are composed of solutions obtained by applying either $1$-interchange \citep{Osman1993} or 2-opt* \citep{Potvin1995} between routes of a given solution. The movements are chosen considering the best improvement strategy. The local search optimizes the solution quality and does not allow infeasible solutions. The feasibility algorithm, on the other hand, attempts to optimize a trade-off between the feasibility gain and the quality of the resulting solution (see \cite{Maximo2021} for more details). If the solution remains infeasible after the neighborhood search, then a new route is added to the solution and the neighborhood search for a feasible solution restarts. This process repeats until either the solution reaches $\overline{m}$ routes or a feasible solution is found. This process does not guarantee the feasibility of an input solution for the HFVRP. 

In both feasibility and local search, the neighborhood investigated for each vertex $v_i$ is limited to the set $ \delta (v_i)$ that contains the $\varphi$ nearest vertices to $v_i$. The parameter $ \varphi \in [1, n-1] $ controls the size of the neighborhood.

\subsection{Control of the Perturbation Degree}
\label{ControlePert} 
The control of the perturbation degree aims to adjust the intensity of the perturbation of a solution. The proposed AILS for the HFVRP considers the same control strategy presented by \cite{Maximo2021}.

Let $d(s,s^r)$ be the distance between the reference solution $s^r$ and the solution $s$ obtained after the local search, namely the perturbation degree of a given perturbation. The employed distance measure  consists of evaluating the number of edges that are different in the two solutions, as presented in Equation \eqref{dab}:

\begin{align}
& d(s,s^r) = |E^s \triangle E^{s^r}|,& \label{dab}
\end{align}
\noindent  where $E^s$ and $ E^{s^r}$ ares the sets of edges of  $s$ and $s^r$, respectively.

The degree of perturbation to be employed by a removal heuristic $\mathscr{H}$ is defined by  $\omega_{\mathscr{H}}$. The greater the value of $\omega_{\mathscr{H}}$, the larger the number of vertices that will be removed in the perturbation step. Therefore, we adjust the $\omega_{\mathscr{H}}$ value to control the diversity of the method. The objective is that the average distance between the solution $s$ and the reference solution $s^r$ represented by $d_{\mathscr{H}}$ to be proportional to the ideal distance $d_{\beta}$, which is a parameter. The larger the $d_{\beta}$, the greater the value of $\omega_{\mathscr{H}}$. The value of $\omega_{\mathscr{H}}$ is adjusted by induction as a function of the real average distance $d_{\mathscr{H}}$ and the ideal distance $d_{\beta}$ as follows $\omega_{\mathscr{H}} = \min\{n,\max\{1, \frac{\omega_{\mathscr{H}}d_{\beta}}{d_{\mathscr{H}}} \}\}$. This adjustment is performed every $\lambda$ choices of heuristic $\mathscr{H}$.

\subsection{Acceptance Criterion}
\label{criterioAceitacao}
The acceptance criteria establish the rules for a given solution $s$ to become the reference solution $s^r$. There are many  possibilities for this step in the literature. In this paper, we calculate a threshold called $\bar{b}$ that provides the minimum quality for a certain solution to replace the reference solution. The threshold value is between the average quality of the solutions obtained after the local search, called $\bar{f}$, and the best solution found during the last $\lambda$ iterations, called $\underline{f}$, i.e. $\bar{b} \in [\bar{f},\underline{f}]$. To calculate  $\bar{f}$, we use the weighted average of the solutions where each new value has a weight of $\lambda^{-1}$, as described in Equation \eqref{fbar}: 


\begin{align}
\bar{f} = \left \{ \begin{array}{ll}
\ \bar{f} (1-\lambda^{-1}) + f(s) \lambda^{-1}, & \textrm{if  $ it > \lambda$}  \\
 (\bar{f} (it-1) + f(s)) it^{-1}, & \textrm{if $ it \leq \lambda$}. 
\end{array}\right. & \label{fbar}
\end{align}

To calculate the value of $\bar{b}$, we use the parameter $\eta \in [0,1]$ which allows the method to control the value of $\bar{b}$ in the interval $[\bar{f},\underline{f}]$: $\bar{b} = \underline{f} + \eta (\bar{f}-\underline{f})$. The greater the value of $\eta$, the more solutions are accepted to become the reference solution. For example, if $\eta=1$, all solutions whose objective function value is higher than the average quality are accepted. If $\eta=0$, only solutions with a quality higher than $\underline{f}$ will be accepted. The acceptance criterion proposed in \citet{Maximo2021} uses the flow of accepted solutions to control the value of $\eta$. In this paper, we choose to use a fixed value of $\eta$.

\section{Computational Experiments}
\label{section:EC}

The AILS heuristic  was coded in Java and we ran computational experiments on a cluster with 104 nodes, each of them with 2 Intel Xeon E5-2680v2 processors with 2.8 GHz, 10 cores and 128 GB DDR3 1866 MHz of RAM. The employed stopping criterion was $40000$ iterations without improvement. To evaluate the performance of AILS we used 12 instances proposed by \cite{Golden1984} and \cite{Taillard1999} whose main characteristics are described in Table \ref{instanciasGolden}. These 12 instances are called Golden in this paper.  AILS was tested in the 12 Golden instances considering every  HFVRP variant, except for HVRPFD and HVRPD variants. Following the same methodology as in the literature, only the four first Golden instances (referred to as `3', `4', `5' and `6') were not considered in the tests of the HVRPFD and HVRPD variants.  Besides the Golden instances, we carried out experiments on the five instances (H1-H5) proposed by \cite{Li2007} and  five instances (N1-N5) introduced by \cite{Brandao2011}, whose characteristics are described in Table \ref{instanciasBrandaoLi}. In these tables, columns A to F refer to the types of vehicles. 

\begin{landscape}
\begin{table}[!ht]
\scalefont{0.55}
\caption{Main characteristics of the instances proposed by \cite{Golden1984} and \cite{Taillard1999}.}
\label{instanciasGolden}
\centering
\begin{tabular}{ccccccccccccccccccccccccccccccc}
\hline
\hline
\multicolumn{2}{c}{Type}& \multicolumn{4}{c}{A} && \multicolumn{4}{c}{B} && \multicolumn{4}{c}{C} && \multicolumn{4}{c}{D} && \multicolumn{4}{c}{E} && \multicolumn{4}{c}{F}\\ 
\cline{3-6} \cline{8-11} \cline{13-16} \cline{18-21} \cline{23-26}\cline{28-31}
Inst.&$n$ & $\bar{q}_A$ & $f_A$ & $r_A$ & $u_A$ && $\bar{q}_B$ & $f_B$ & $r_B$ & $u_B$ && $\bar{q}_C$ & $f_C$ & $r_C$ & $u_C$ && $\bar{q}_D$ & $f_D$ & $r_D$ & $u_D$ && $\bar{q}_E$ & $f_E$ & $r_E$ & $u_E$ && $\bar{q}_F$ & $f_F$ & $r_F$ & $u_F$ \\
\hline
\rowcolor{Gray}
3&20&20&20&1.0&20&&30&35&1.1&20&&40&50&1.2&20&&70&120&1.7&20&&120&225&2.5&20&&&&&\\
4&20&60&1000&1.0&20&&80&1500&1.1&20&&150&3000&1.4&20&&&&&&&&&&&&&&&\\
\rowcolor{Gray}
5&20&20&20&1.0&20&&30&35&1.1&20&&40&50&1.2&20&&70&120&1.7&20&&120&225&2.5&20&&&&&\\
6&20&60&1000&1.0&20&&30&1500&1.1&20&&150&3000&1.4&20&&&&&&&&&&&&&&&\\
\rowcolor{Gray}
13&50&20&20&1.0&4&&30&35&1.1&2&&40&50&1.2&4&&70&120&1.7&4&&120&225&2.5&2&&200&400&3.2&1\\
14&50&120&1000&1.0&4&&160&1500&1.1&2&&300&3500&1.4&1&&&&&&&&&&&&&&&\\
\rowcolor{Gray}
15&50&50&100&1.0&4&&100&250&1.6&3&&160&450&2.0&2&&&&&&&&&&&&&&&\\
16&50&40&100&1.0&2&&80&200&1.6&4&&140&400&2.1&3&&&&&&&&&&&&&&&\\
\rowcolor{Gray}
17&75&50&25&1.0&4&&120&80&1.2&4&&200&150&1.5&2&&350&320&1.8&1&&&&&&&&&&\\
18&75&20&10&1.0&4&&50&35&1.3&4&&100&100&1.9&2&&150&180&2.4&2&&250&400&2.9&1&&400&800&3.2&1\\
\rowcolor{Gray}
19&100&100&500&1.0&4&&200&1200&1.4&3&&300&2100&1.7&3&&&&&&&&&&&&&&&\\
20&100&60&100&1.0&6&&140&300&1.7&4&&200&500&2.0&3&&&&&&&&&&&&&&&\\
\hline
\hline
\end{tabular}
\end{table}

\begin{table}[!ht]
\scalefont{0.55}
\caption{Main characteristics of the instances proposed by \cite{Li2007} and \cite{Brandao2011}.}
\label{instanciasBrandaoLi}
\centering
\begin{tabular}{ccccccccccccccccccccccccc}
\hline
\hline
\multicolumn{2}{c}{Type}& \multicolumn{3}{c}{A} && \multicolumn{3}{c}{B} && \multicolumn{3}{c}{C} && \multicolumn{3}{c}{D} && \multicolumn{3}{c}{E} && \multicolumn{3}{c}{F}\\ 
\cline{3-5} \cline{7-9} \cline{11-13} \cline{15-17} \cline{19-21}\cline{23-25}
Inst.&$n$ & $\bar{q}_A$ & $r_A$ & $u_A$ && $\bar{q}_B$ & $r_B$ & $u_B$ && $\bar{q}_C$ & $r_C$ & $u_C$ && $\bar{q}_D$ & $r_D$ & $u_D$ && $\bar{q}_E$ & $r_E$ & $u_E$ && $\bar{q}_F$ & $r_F$ & $u_F$ \\
\hline
\rowcolor{Gray}
N1&150&50&1&5&&100&1.5&4&&150&1.9&4&&200&2.2&3&&250&2.6&2&&&&\\
N2&199&50&1&8&&100&1.5&6&&150&1.9&5&&200&2.2&4&&250&2.6&2&&350&3.2&1\\
\rowcolor{Gray}
N3&120&50&1&6&&100&1.5&3&&150&1.9&3&&200&2.2&2&&&&&&&&\\
N4&100&50&1&4&&120&1.6&4&&180&2.1&4&&240&2.6&2&&&&&&&&\\
\rowcolor{Gray}
N5&134&900&1&5&&1500&1.5&3&&2000&1.8&2&&2500&2.2&1&&&&&&&&\\
H1&200&50&1&8&&100&1.1&6&&200&1.2&4&&500&1.7&3&&1000&2.5&1&&&&\\
\rowcolor{Gray}
H2&240&50&1&10&&100&1.1&5&&200&1.2&5&&500&1.7&4&&1000&2.5&1&&&&\\
H3&280&50&1&10&&100&1.1&5&&200&1.2&5&&500&1.7&4&&1000&2.5&2&&&&\\
\rowcolor{Gray}
H4&320&50&1&10&&100&1.1&8&&200&1.2&5&&500&1.7&2&&1000&2.5&2&&1500&3&1\\
H5&360&50&1&10&&100&1.2&8&&200&1.5&5&&500&1.8&1&&1500&2.5&2&&2000&3&1\\
\hline
\hline
\end{tabular}
\end{table}

\end{landscape}

\subsection{Parameter tuning}\label{parametertuning}
 In summary, AILS has five parameters:

\begin{itemize}
    \item $\alpha$: parameter required by the perturbation heuristics that indicates the probability of varying the number of routes and vehicles assigned to each route;
    

    \item $d_{\beta}$: reference distance between the reference solution and the current solution obtained after the local search; 
 
    \item $\eta$: hyperparameter to define the threshold $\bar{b}$ in the acceptance criterion;
  
    \item $\gamma$: number of iterations for AILS to perform a new adjustment of parameter $\omega$;

    \item $\varphi$: parameter of the feasibility and local search that refers to the maximum cardinality of $\delta(v)$ -- nearest neighbors of $v$.

\end{itemize}

The parameter tuning of AILS considered the Golden instances applied to all five HFVRP variants.  The value analyzed as a reference for parameter adjustment was the average gap of 10 executions ($avg$). The mean gap represents the relative distance of the quality of the solutions found in relation to the best-known solutions (BKS) and is calculated according to Equation \eqref{gap}:

\begin{align}
&gap=100(avg-BKS)/BKS.& \label{gap}
\end{align}

The values considered for each parameter are

\begin{itemize}
    \item $\alpha \in \{ 0.1, 0.2,\ldots, 1\}$
     \item $d_{\beta} \in \{5,10,\ldots, 40\}$
    \item $\eta \in \{0.1, 0.2, \ldots, 0.6\}$
    \item  $\gamma \in \{10,20,\ldots,100\}$
    \item $\varphi \in \{10,20,\ldots,100\}$.
\end{itemize}

The parameterization method iteratively varied each of the parameters in their respective interval and fixed the values of the other parameters. The priority order to set up the parameter values was: $\alpha, d_{\beta}, \eta, \gamma$ and $\varphi$. Initially, the parameters were set with the following configuration: $\alpha=0.4$, $d_{\beta}=15$, $\eta=0.2$, $\gamma=30$ and $ \varphi=60$. These values were obtained in preliminary experiments by optimizing the gaps. 
Table \ref{parametros} presents the best configuration according to this parameter tuning methodology.
\begin{table}[!ht]
\caption{Values employed in the AILS after the set up of parameters.}
\label{parametros}
\centering
\begin{tabular}{ccccc}
\hline
\hline
Parameter & Type & Interval & Value \\
\rowcolor{Gray}
$\alpha$ & Real & [0.1, 1] & 0.4\\
$d_{\beta}$ & Integer & [5, 50] & 15\\
\rowcolor{Gray}
$\eta$ & Real & [0.1, 1] & 0.2\\
$\gamma$ & Integer & [10, 100] & 20\\
\rowcolor{Gray}
$\varphi$ & Integer & [10, 100] & 20\\
\hline
\hline
\end{tabular}
\end{table}

\begin{figure}[h!t]
\center
\begin{tikzpicture}
\pgfplotsset{every axis legend/.append style={
at={(0.5,1.05)},
anchor=south},
}
\begin{axis}[
width=6cm,
height=4.5cm,
legend columns=2,
xmin=0.05,
xmax=1.05,
xtick={0.1, 0.3, 0.5, 0.7,0.9 },
ymin=0,
ymax=0.28,
grid=major,
xlabel={$\alpha$},
ylabel={$gap$}
]
\addplot+[
red!70!white,
mark options={fill=red!70!white},
line width=2pt,
mark=triangle*,
mark size=1.5pt
]
table[x=x,y=y] {dados/Tunning/variandoAlfa.txt};
\legend{}
\end{axis}
\end{tikzpicture}
\begin{tikzpicture}
\pgfplotsset{every axis legend/.append style={
at={(0.5,1.05)},
anchor=south},
}
\begin{axis}[
width=6cm,
height=4.5cm,
legend columns=2,
xmin=3,
xmax=42,
xtick={5,15, 25, 35},
ymin=0,
ymax=0.28,
grid=major,
xlabel={$d_{\beta}$},
ylabel={}
]
\addplot+[
blue!50!white,
mark options={fill=blue!50!white},
line width=2pt,
mark=*,
mark size=1.5pt
]
table[x=x,y=y] {dados/Tunning/variandoDBeta.txt};
\legend{}
\end{axis}
\end{tikzpicture}
\begin{tikzpicture}
\pgfplotsset{every axis legend/.append style={
at={(0.5,1.05)},
anchor=south},
}
\begin{axis}[
width=6cm,
height=4.5cm,
legend columns=2,
xmin=0.05,
xmax=0.65,
xtick={0.1, 0.2, 0.3,0.4, 0.5,0.6},
ymin=0,
ymax=0.28,
grid=major,
xlabel={$\eta$},
ylabel={}
]
\addplot+[
violet!50!white,
mark options={fill=violet!50!white},
line width=2pt,
mark=diamond*,
mark size=1.5pt
]
table[x=x,y=y] {dados/Tunning/variandoEta.txt};
\legend{}
\end{axis}
\end{tikzpicture}\\
\begin{tikzpicture}
\pgfplotsset{every axis legend/.append style={
at={(0.5,1.05)},
anchor=south},
}
\begin{axis}[
width=6cm,
height=4.5cm,
legend columns=2,
xmin=5,
xmax=105,
xtick={10, 30,50,70,90},
ymin=0,
ymax=0.28,
grid=major,
xlabel={$\gamma$},
ylabel={}
]
\addplot+[
orange!70!white,
mark options={fill=orange!70!white},
line width=2pt,
mark=square*,
mark size=1.5pt
]
table[x=x,y=y] {dados/Tunning/variandoGamma.txt};
\legend{}
\end{axis}
\end{tikzpicture}
\begin{tikzpicture}
\pgfplotsset{every axis legend/.append style={
at={(0.5,1.05)},
anchor=south},
}
\begin{axis}[
width=6cm,
height=4.5cm,
legend columns=2,
xmin=5,
xmax=105,
xtick={10, 30, 50, 70, 90},
ymin=0,
ymax=0.28,
grid=major,
xlabel={$\varphi$},
ylabel={}
]
\addplot+[
green!50!white,
mark options={fill=green!50!white},
line width=2pt,
mark=square*,
mark size=1.5pt
]
table[x=x,y=y] {dados/Tunning/variandoVarphi.txt};
\legend{}
\end{axis}
\end{tikzpicture}
\caption{Results achieved in the fine tuning of the AILS parameters.}
\label{fig:variandoParam}
\end{figure}
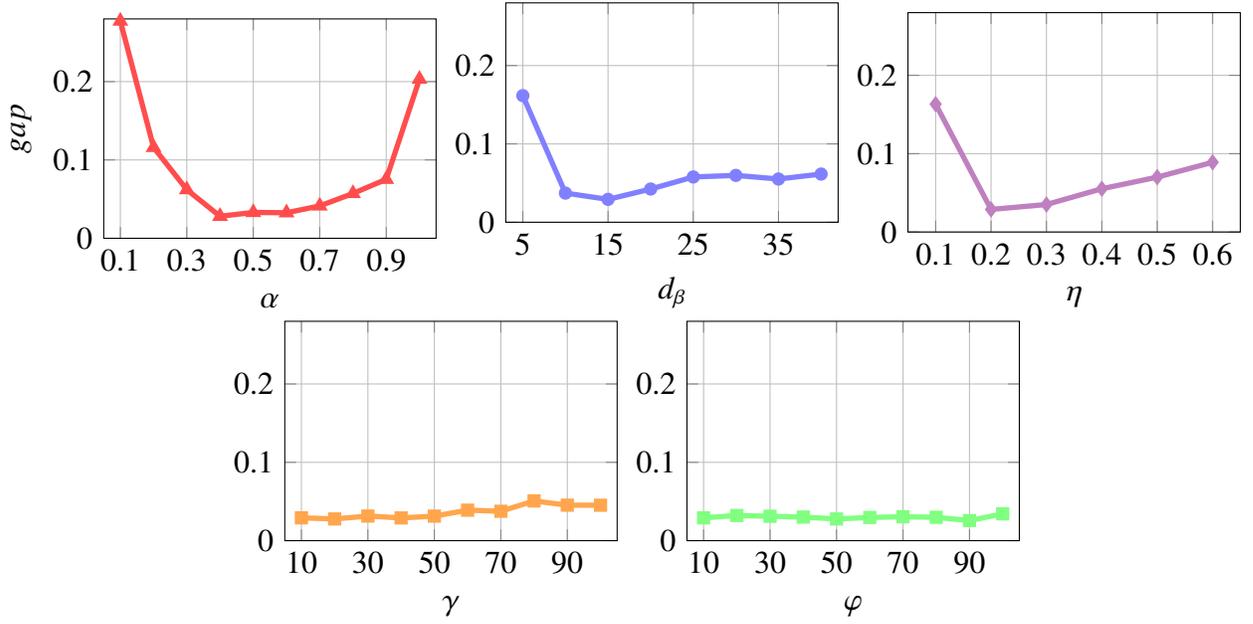

Figure \ref{fig:variandoParam} illustrates the  parameter setting. It shows that $\alpha$, $\eta$ and $d_{\beta}$ have greater sensitivity in comparison to the other parameters.



\subsection{Experiments}
In these experiments, we compared the AILS results with other algorithms in the literature: ILS-RVND \citep{Penna2013}, ILS-RVND-SP \citep{Subramanian2012} and UHGS \citep{Vidal2014}. Table \ref{comparacao} summarizes the information about the setup of the experiments carried out by the authors: the programming language of the metaheuristic;  the clock rate of the processor on which the authors ran their experiments; the number of runs for each tested instance; and the variants the authors considered in their experiments by the proposed metaheuristics. This table also shows the same information regarding the experiments we carried out with AILS.


\begin{table}[!ht]
\caption{Experimental setup of the analyzed algorithms.}
\label{comparacao}
\centering
\begin{tabular}{lcccc}
\hline
\hline
\rowcolor{Gray} 
Method   & ILS-RVND  & ILS-RVND-SP   & UHGS  & AILS\\
\hline
Language   & C++       & C++           & C++   & JAVA\\
\rowcolor{Gray}
Processor & 2.93 GHz & 2.93 GHz & 2.4 GHz &  2.8 GHz \\
Nb. of runs     & 30 &10 & 10 & 10\\
\rowcolor{Gray}
HVRPFD      & \checkmark&\checkmark&&\checkmark\\
HVRPD       & \checkmark&\checkmark&&\checkmark\\
\rowcolor{Gray}
FSMFD       &\checkmark&\checkmark&\checkmark&\checkmark\\
FSMF        &\checkmark&\checkmark&\checkmark&\checkmark\\
\rowcolor{Gray}
FSMD        &\checkmark&\checkmark&\checkmark&\checkmark\\
\hline
\hline
\end{tabular}
\end{table}

Tables \ref{GoldenHVRPFD} to \ref{BrandaoFSMD} show the results of the experiments. Besides `BKS', `$avg$' and `$gap$', these tables report the best solutions in all runs (`Best') and the average time in seconds of the algorithms (`Time').

    
    

    

Next, we analyze the results for each of the variants.

\begin{itemize}
    \item \textbf{HVRPFD}: Table~\ref{GoldenHVRPFD} shows the results of AILS,  ILS-RVND and ILS-RVND-SP considering HVRPFD on the Golden instances. According to the presented results, AILS achieved the best average gap in all instances. The best AILS solution was equal to the best-known solution (BKS) in all instances. ILS-RVND-SP was faster than the other heuristics on all instances. ILS-RVND-SP presented the second best gap in seven of the eight instances. 
    
    \item \textbf{HVRPD}: Table \ref{GoldenHVRPD} shows the results of AILS, ILS-RVND and ILS-RVND-SP considering HVRPD on the Golden instances. The results for this variant were similar to those presented for the HVRPFD. The average AILS gap was better on all instances. AILS also found the best solution in all cases. The time to obtain the solutions was better in only one instance. ILS-RVND-SP achieved the best average gap on two instances and the same best solutions as AILS. ILS-RVND-SP was the fastest method on seven instances. ILS-RVND had the best average gap on one instance and the longest time on six instances. This variant was also analyzed using the instances proposed by \cite{Brandao2011} and \cite{Li2007}. Table~\ref{BrandaoLiHVRPFD} shows the results of ILS-RVND-SP and AILS. It is worth pointing out that AILS obtained a better solution than those already known for six instances out of the 10 considered, therefore, updating the literature BKS of six instances. AILS had the tightest average gap in all 10 instances and the longest time in four instances. 
    
    \item \textbf{FSMFD}:  Table \ref{GoldenFSMFD} shows the results of AILS with ILS-RVND, ILS-RVND-SP and UHGS considering FSMFD on the Golden instances. The average AILS gap was the tightest among the compared methods on seven of the 12 instances. The best solutions obtained by AILS were equal to the BKS on ten instances. AILS presented the best time on four instances and the second shortest time on four instances. ILS-RVND-SP outperformed AILS on four instances. ILS-RVND-SP also had the best time on eight instances and the  best solutions on all instances. UHGS had a null average gap on eight instances, however, the average time to obtain the solution was significantly higher than the other three algorithms.

    \item \textbf{FSMF}: Table \ref{GoldenFSMF} presents the results of AILS with ILS-RVND, ILS-RVND-SP and UHGS considering FSMF on the Golden instances. AILS had the best average result on eight of the 12 instances. The average gap of the solutions obtained by AILS was $0.0297$, whereas the gaps of UHGS, ILS-RVND-SP and ILS-RVND were, respectively, $0.0396$, $0.1244$ and $0.2527$. The best solutions found by AILS were the same as the BKSs on all tested instances. UHGS had the best average gap on seven instances and the second best gap on four instances. ILS-RVND-SP presented very good solutions considering it took, on average, half the time of AILS. 
    
    \item \textbf{FSMD}: Table \ref{GoldenFSMD} shows the results of AILS with ILS-RVND, ILS-RVND-SP and UHGS considering FSMD on the Golden instances. AILS presented a null gap on 8 of the 12 instances and the best results on 9 of them. UHGS presented a null gap on instance `16', whereas AILS obtained an average solution whose gap was $0.1627$. AILS had the lowest times on three of the 12 instances and found the BKS on 11 of them. ILS-RVND-SP achieved larger gaps than AILS, however, its average time was much lower. This variant was also analyzed using the instances proposed by \cite{Brandao2011}. Table \ref{BrandaoFSMD} presents the results of AILS and ILS-RVND-SP. AILS was able to find a better BKS on the `N1' and `N2' instances. These solutions were obtained in preliminary experiments. On the five instances, AILS had the best average results. ILS-RVND-SP presented a better best solution than AILS on instance `N1'.

\end{itemize}

\begin{table}[!ht]
\centering
\scalefont{0.55}
\caption{Results of the analyzed methods for HVRPFD on the Golden instances.}
\label{GoldenHVRPFD}
\centering
\begin{tabular}{ccccccccccccc}
\hline
&& \multicolumn{3}{c}{ILS-RVND \citep{Penna2013}} && \multicolumn{3}{c}{ILS-RVND-SP \citep{Subramanian2012}} && \multicolumn{3}{c}{AILS}\\ 
\cline{3-5} \cline{7-9} \cline{11-13} 
Inst.&BKS&Avg~(\textit{gap})&Best&Time&&Avg~(\textit{gap})&Best&Time&&Avg~(\textit{gap})&Best&Time\\
\hline
\rowcolor{Gray}
13&$\mathbf{3185.09^*}$&3189.17~(0.1281)&\textbf{3185.09}&19.04&&3186.32~(0.0386)&\textbf{3185.09}&\textbf{1.99}&&\textbf{3185.09~(0.0000)}&\textbf{3185.09}&6.92\\
14&$\mathbf{10107.53^*}$&10107.94~(0.0041)&\textbf{10107.53}&11.28&&10110.61~(0.0305)&\textbf{10107.53}&\textbf{1.29}&&\textbf{10107.53~(0.0000)}&\textbf{10107.53}&2.98\\
\rowcolor{Gray}
15&$\mathbf{3065.29^*}$&3065.34~(0.0016)&\textbf{3065.29}&12.48&&\textbf{3065.29~(0.0000)}&\textbf{3065.29}&\textbf{1.77}&&\textbf{3065.29~(0.0000)}&\textbf{3065.29}&26.70\\
16&$\mathbf{3265.41^*}$&3278.06~(0.3874)&\textbf{3265.41}&12.22&&3273.15~(0.2370)&\textbf{3265.41}&\textbf{1.67}&&\textbf{3265.41~(0.0000)}&\textbf{3265.41}&25.41\\
\rowcolor{Gray}
17&$\mathbf{2076.96^*}$&2083.19~(0.3000)&\textbf{2076.96}&29.59&&2081.55~(0.2210)&\textbf{2076.96}&\textbf{5.95}&&\textbf{2076.96~(0.0000)}&\textbf{2076.96}&19.63\\
18&$\mathbf{3743.58^*}$&3758.84~(0.4076)&\textbf{3743.58}&36.38&&3758.83~(0.4074)&\textbf{3743.58}&\textbf{16.47}&&\textbf{3743.58~(0.0000)}&\textbf{3743.58}&35.38\\
\rowcolor{Gray}
19&$10420.34$&10421.39~(0.0101)&\textbf{10420.34}&73.66&&10421.05~(0.0068)&\textbf{10420.34}&\textbf{15.8}&&\textbf{10420.34~(0.0000)}&\textbf{10420.34}&41.23\\
20&$4760.68$&4839.53~(1.6563)&4788.49&68.46&&4822.16~(1.2914)&4761.26&\textbf{16.87}&&\textbf{4764.34~(0.0769)}&\textbf{4760.68}&60.27\\
\rowcolor{Gray}
Avg&&~(0.3619)&&32.89&&~(0.2791)&&7.73&&~(0.0096)&&27.32\\
\hline
\end{tabular}
\end{table}

\begin{table}[!ht]
\centering
\scalefont{0.55}
\caption{Results of the analyzed methods for HVRPD on the Golden instances.}
\label{GoldenHVRPD}
\centering
\begin{tabular}{ccccccccccccc}
\hline
&& \multicolumn{3}{c}{ILS-RVND \citep{Penna2013}} && \multicolumn{3}{c}{ILS-RVND-SP \citep{Subramanian2012}} && \multicolumn{3}{c}{AILS}\\ 
\cline{3-5} \cline{7-9} \cline{11-13} 
Inst.&BKS&Avg~(\textit{gap})&Best&Time&&Avg~(\textit{gap})&Best&Time&&Avg~(\textit{gap})&Best&Time\\
\hline
\rowcolor{Gray}
13&$\mathbf{1517.84^*}$&1518.58~(0.0488)&\textbf{1517.84}&19.29&&\textbf{1517.84~(0.0000)}&\textbf{1517.84}&\textbf{1.33}&&\textbf{1517.84~(0.0000)}&\textbf{1517.84}&3.93\\
14&$\mathbf{607.53^*}$&607.64~(0.0181)&\textbf{607.53}&11.20&&608.74~(0.1992)&\textbf{607.53}&\textbf{1.09}&&\textbf{607.53~(0.0000)}&\textbf{607.53}&2.62\\
\rowcolor{Gray}
15&$\mathbf{1015.29^*}$&1015.33~(0.0039)&\textbf{1015.29}&12.56&&\textbf{1015.29~(0.0000)}&\textbf{1015.29}&\textbf{2.13}&&\textbf{1015.29~(0.0000)}&\textbf{1015.29}&22.90\\
16&$\mathbf{1144.94^*}$&1145.04~(0.0087)&\textbf{1144.94}&12.29&&1145.23~(0.0253)&\textbf{1144.94}&\textbf{1.41}&&\textbf{1144.94~(0.0000)}&\textbf{1144.94}&1.86\\
\rowcolor{Gray}
17&$\mathbf{1061.96^*}$&1065.27~(0.3117)&\textbf{1061.96}&29.92&&1064.67~(0.2552)&\textbf{1061.96}&\textbf{4.22}&&\textbf{1062.17~(0.0198)}&\textbf{1061.96}&23.71\\
18&$\mathbf{1823.58^*}$&1832.52~(0.4902)&\textbf{1823.58}&38.34&&1831.48~(0.4332)&\textbf{1823.58}&\textbf{4.06}&&\textbf{1824.86~(0.0702)}&\textbf{1823.58}&62.82\\
\rowcolor{Gray}
19&$1117.51$&\textbf{1120.34~(0.2532)}&\textbf{1120.34}&67.72&&1121.11~(0.3221)&\textbf{1120.34}&\textbf{9.12}&&\textbf{1120.34~(0.2532)}&\textbf{1120.34}&46.06\\
20&$\mathbf{1534.17^*}$&1544.08~(0.6460)&\textbf{1534.17}&63.77&&1536.89~(0.1773)&\textbf{1534.17}&8.89&&\textbf{1534.17~(0.0000)}&\textbf{1534.17}&\textbf{4.14}\\
\rowcolor{Gray}
Avg&&~(0.2226)&&31.89&&~(0.1765)&&4.03&&~(0.0429)&&21.01\\
\hline
\end{tabular}
\end{table}

\begin{table}[!ht]
\centering
\scalefont{0.55}
\caption{Results of the analyzed methods for HVRPD on the instances introduced by \cite{Brandao2011} and \cite{Li2007}.}
\label{BrandaoLiHVRPFD}
\centering
\begin{threeparttable}
\begin{tabular}{ccccccccc}
\hline
&& \multicolumn{3}{c}{ILS-RVND-SP \citep{Subramanian2012}} && \multicolumn{3}{c}{AILS}\\ 
\cline{3-5} \cline{7-9} 
Inst.&BKS&Avg~(\textit{gap})&Best&Time&&Avg~(\textit{gap})&Best&Time\\
\hline
\rowcolor{Gray}
N1&$2233.90$\tnote{1}&2244.31~(0.4660)&2235.87&\textbf{51.50}&&\textbf{2237.43~(0.1580)}&\textbf{2235.34}&61.41\\
N2&$2851.94$\tnote{1}&2906.24~(1.9040)&2864.83&102.77&&\textbf{2861.97~(0.3517)}&\textbf{2852.81}&\textbf{73.71}\\
\rowcolor{Gray}
N3&$2378.99$\tnote{2}&2382.1~(0.1307)&\textbf{2378.99}&51.71&&\textbf{2379.99~(0.0420)}&\textbf{2378.99}&\textbf{27.87}\\
N4&$1839.22$&\textbf{1839.22~(0.0000)}&\textbf{1839.22}&9.64&&\textbf{1839.22~(0.0000)}&\textbf{1839.22}&\textbf{7.20}\\
\rowcolor{Gray}
N5&$2047.81$\tnote{2}&\textbf{2047.81~(0.0000)}&\textbf{2047.81}&52.33&&\textbf{2047.81~(0.0000)}&\textbf{2047.81}&\textbf{19.30}\\
H1&$12050.08$&12052.69~(0.0217)&\textbf{12050.08}&\textbf{72.10}&&\textbf{12052.08~(0.0166)}&\textbf{12050.08}&82.12\\
\rowcolor{Gray}
H2&$10130.30$\tnote{1}&10436.2~(3.0197)&10329.15&\textbf{176.43}&&\textbf{10198.87~(0.6769)}&\textbf{10161.90}&203.79\\
H3&$16192.26$\tnote{1}&16526.89~(2.0666)&16282.41&259.61&&\textbf{16205.86~(0.0840)}&\textbf{16193.48}&\textbf{188.69}\\
\rowcolor{Gray}
H4&$17273.75$\tnote{1}&18022.37~(4.3339)&17743.68&384.52&&\textbf{17459.21~(1.0737)}&\textbf{17397.28}&\textbf{228.16}\\
H5&$23024.58$\tnote{1}&23948.97~(4.0148)&23493.87&\textbf{621.17}&&\textbf{23222.35~(0.8590)}&\textbf{23088.98}&703.71\\
\rowcolor{Gray}
Avg&&~(1.5957)&&178.18&&~(0.3262)&&159.60\\
\hline
\end{tabular}
\begin{tablenotes}
\item[1] Best solutions found by AILS.
\item[2] Best solutions found by ILS-RVND-SP.
\end{tablenotes}
\end{threeparttable}
\end{table}

\begin{landscape}
\begin{table}[!ht]
\centering
\scalefont{0.6}
\caption{Results of the analyzed methods for FSMFD on the Golden instances.}
\label{GoldenFSMFD}
\centering
\begin{tabular}{crccccccccccccccc}
\hline
  && \multicolumn{3}{c}{ILS-RVND \citep{Penna2013}} &&  \multicolumn{3}{c}{ILS-RVND-SP \citep{Subramanian2012}} &  \multicolumn{3}{c}{UHGS \citep{Vidal2014}} && \multicolumn{3}{c}{AILS}\\ 
\cline{3-5} \cline{7-9} \cline{11-13} \cline{15-17} 
Inst.&BKS&Avg~(\textit{gap})&Best&Time&&Avg~(\textit{gap})&Best&Time&&Avg~(\textit{gap})&Best&Time&&Avg~(\textit{gap})&Best&Time\\
\hline
\rowcolor{Gray}
3&$\mathbf{1144.22^*}$&\textbf{1144.22~(0.0000)}&\textbf{1144.22}&4.05&&\textbf{1144.22~(0.0000)}&\textbf{1144.22}&\textbf{0.34}&&\textbf{1144.22~(0.0000)}&\textbf{1144.22}&10.20&&\textbf{1144.22~(0.0000)}&\textbf{1144.22}&9.86\\
4&$\mathbf{6437.33^*}$&6437.66~(0.0051)&\textbf{6437.33}&3.03&&\textbf{6437.33~(0.0000)}&\textbf{6437.33}&\textbf{0.31}&&\textbf{6437.33~(0.0000)}&\textbf{6437.33}&13.80&&\textbf{6437.33~(0.0000)}&\textbf{6437.33}&8.76\\
\rowcolor{Gray}
5&$\mathbf{1322.26^*}$&\textbf{1322.26~(0.0000)}&\textbf{1322.26}&4.85&&\textbf{1322.26~(0.0000)}&\textbf{1322.26}&0.28&&\textbf{1322.26~(0.0000)}&\textbf{1322.26}&10.20&&\textbf{1322.26~(0.0000)}&\textbf{1322.26}&\textbf{0.14}\\
6&$\mathbf{6516.47^*}$&\textbf{6516.47~(0.0000)}&\textbf{6516.47}&3.01&&\textbf{6516.47~(0.0000)}&\textbf{6516.47}&0.32&&\textbf{6516.47~(0.0000)}&\textbf{6516.47}&13.80&&\textbf{6516.47~(0.0000)}&\textbf{6516.47}&\textbf{0.10}\\
\rowcolor{Gray}
13&$\mathbf{2964.65^*}$&2971.32~(0.2250)&\textbf{2964.65}&27.44&&\textbf{2964.65~(0.0000)}&\textbf{2964.65}&\textbf{1.70}&&\textbf{2964.65~(0.0000)}&\textbf{2964.65}&30.60&&2964.83~(0.0061)&\textbf{2964.65}&16.60\\
14&$\mathbf{9126.90^*}$&9126.91~(0.0001)&\textbf{9126.90}&11.66&&\textbf{9126.90~(0.0000)}&\textbf{9126.90}&1.53&&\textbf{9126.90~(0.0000)}&\textbf{9126.90}&47.40&&\textbf{9126.90~(0.0000)}&\textbf{9126.90}&\textbf{0.97}\\
\rowcolor{Gray}
15&$\mathbf{2634.96^*}$&2635.02~(0.0023)&\textbf{2634.96}&13.83&&\textbf{2634.96~(0.0000)}&\textbf{2634.96}&\textbf{1.34}&&2635.06~(0.0038)&\textbf{2634.96}&42.60&&\textbf{2634.96~(0.0000)}&\textbf{2634.96}&21.40\\
16&$\mathbf{3168.92^*}$&3170.81~(0.0596)&\textbf{3168.92}&18.20&&\textbf{3168.92~(0.0000)}&\textbf{3168.92}&\textbf{6.72}&&\textbf{3168.92~(0.0000)}&\textbf{3168.92}&48.00&&3168.94~(0.0006)&\textbf{3168.92}&25.52\\
\rowcolor{Gray}
17&$\mathbf{2004.48^*}$&2012.23~(0.3866)&\textbf{2004.48}&43.68&&2007.12~(0.1317)&\textbf{2004.48}&\textbf{6.96}&&2007.04~(0.1277)&\textbf{2004.48}&79.80&&\textbf{2004.48~(0.0000)}&\textbf{2004.48}&24.72\\
18&$\mathbf{3147.99^*}$&3158.24~(0.3256)&3149.63&47.80&&\textbf{3148.91~(0.0292)}&\textbf{3147.99}&\textbf{4.21}&&3148.99~(0.0318)&3148.99&76.80&&3151.40~(0.1083)&3148.56&45.59\\
\rowcolor{Gray}
19&$\mathbf{8661.81^*}$&8664.81~(0.0346)&\textbf{8661.81}&59.13&&\textbf{8662.89~(0.0125)}&\textbf{8661.81}&\textbf{29.86}&&8663.04~(0.0142)&\textbf{8661.81}&234.60&&8663.26~(0.0167)&8661.82&40.72\\
20&$4153.02$&4155.90~(0.0693)&\textbf{4153.02}&59.07&&4153.12~(0.0024)&\textbf{4153.02}&37.21&&\textbf{4153.02~(0.0000)}&\textbf{4153.02}&103.80&&4153.12~(0.0024)&\textbf{4153.02}&\textbf{25.42}\\
\rowcolor{Gray}
Avg&&~(0.0924)&&24.65&&~(0.0147)&&7.57&&~(0.0148)&&59.30&&~(0.0112)&&18.32\\
\hline
\end{tabular}
\end{table}
\end{landscape}

\begin{landscape}
\begin{table}[!ht]
\centering
\scalefont{0.6}
\caption{Results of the analyzed methods for FSMF on the Golden instances.}
\label{GoldenFSMF}
\centering
\begin{tabular}{crccccccccccccccc}
\hline
  && \multicolumn{3}{c}{ILS-RVND \citep{Penna2013}} &&  \multicolumn{3}{c}{ILS-RVND-SP \citep{Subramanian2012}} &  \multicolumn{3}{c}{UHGS \citep{Vidal2014}} && \multicolumn{3}{c}{AILS}\\ 
\cline{3-5} \cline{7-9} \cline{11-13} \cline{15-17} 
Inst.&BKS&Avg~(\textit{gap})&Best&Time&&Avg~(\textit{gap})&Best&Time&&Avg~(\textit{gap})&Best&Time&&Avg~(\textit{gap})&Best&Time\\
\hline
\rowcolor{Gray}
3&$\mathbf{961.03^*}$&961.10~(0.0073)&\textbf{961.03}&4.91&&\textbf{961.03~(0.0000)}&\textbf{961.03}&\textbf{0.28}&&\textbf{961.03~(0.0000)}&\textbf{961.03}&12.00&&\textbf{961.03~(0.0000)}&\textbf{961.03}&0.41\\
4&$\mathbf{6437.33^*}$&6437.63~(0.0047)&\textbf{6437.33}&3.16&&\textbf{6437.33~(0.0000)}&\textbf{6437.33}&\textbf{0.25}&&\textbf{6437.33~(0.0000)}&\textbf{6437.33}&13.80&&\textbf{6437.33~(0.0000)}&\textbf{6437.33}&7.13\\
\rowcolor{Gray}
5&$\mathbf{1007.05^*}$&\textbf{1007.05~(0.0000)}&\textbf{1007.05}&5.88&&1008.76~(0.1698)&\textbf{1007.05}&\textbf{0.25}&&\textbf{1007.05~(0.0000)}&\textbf{1007.05}&13.80&&1007.20~(0.0149)&\textbf{1007.05}&11.70\\
6&$\mathbf{6516.47^*}$&\textbf{6516.47~(0.0000)}&\textbf{6516.47}&3.07&&\textbf{6516.47~(0.0000)}&\textbf{6516.47}&\textbf{0.20}&&\textbf{6516.47~(0.0000)}&\textbf{6516.47}&13.80&&\textbf{6516.47~(0.0000)}&\textbf{6516.47}&0.70\\
\rowcolor{Gray}
13&$\mathbf{2406.36^*}$&2419.38~(0.5411)&2408.41&30.29&&2411.31~(0.2057)&\textbf{2406.36}&\textbf{1.96}&&2406.57~(0.0087)&\textbf{2406.36}&61.20&&\textbf{2406.36~(0.0000)}&\textbf{2406.36}&27.31\\
14&$\mathbf{9119.03^*}$&\textbf{9119.03~(0.0000)}&\textbf{9119.03}&11.89&&\textbf{9119.03~(0.0000)}&\textbf{9119.03}&\textbf{1.64}&&\textbf{9119.03~(0.0000)}&\textbf{9119.03}&52.80&&\textbf{9119.03~(0.0000)}&\textbf{9119.03}&24.84\\
\rowcolor{Gray}
15&$\mathbf{2586.37^*}$&2586.80~(0.0166)&\textbf{2586.37}&20.24&&\textbf{2586.37~(0.0000)}&\textbf{2586.37}&\textbf{6.02}&&\textbf{2586.37~(0.0000)}&\textbf{2586.37}&43.80&&2586.42~(0.0019)&\textbf{2586.37}&8.64\\
16&$\mathbf{2720.43^*}$&2737.59~(0.6308)&\textbf{2720.43}&20.67&&2724.55~(0.1514)&\textbf{2720.43}&\textbf{3.85}&&\textbf{2720.43~(0.0000)}&\textbf{2720.43}&39.60&&\textbf{2720.43~(0.0000)}&\textbf{2720.43}&18.77\\
\rowcolor{Gray}
17&$\mathbf{1734.53^*}$&1748.06~(0.7800)&\textbf{1734.53}&52.49&&1744.23~(0.5592)&\textbf{1734.53}&\textbf{11.61}&&1735.37~(0.0484)&\textbf{1734.53}&105.00&&\textbf{1734.69~(0.0092)}&\textbf{1734.53}&38.18\\
18&$\mathbf{2369.65^*}$&2380.98~(0.4781)&2371.48&55.35&&\textbf{2373.79~(0.1747)}&\textbf{2369.65}&\textbf{11.83}&&2374.16~(0.1903)&\textbf{2369.65}&103.80&&2375.95~(0.2659)&\textbf{2369.65}&32.62\\
\rowcolor{Gray}
19&$8661.81$&8665.31~(0.0404)&8662.86&63.92&&\textbf{8662.54~(0.0084)}&\textbf{8661.81}&\textbf{25.15}&&8663.97~(0.0249)&8662.86&222.00&&8662.73~(0.0106)&\textbf{8661.81}&35.48\\
20&$\mathbf{4029.61^*}$&4051.11~(0.5336)&4037.90&93.88&&4038.63~(0.2238)&4032.81&46.06&&4037.77~(0.2025)&4034.42&135.60&&\textbf{4031.78~(0.0539)}&\textbf{4029.61}&\textbf{34.08}\\
\rowcolor{Gray}
Avg&&~(0.2527)&&30.48&&~(0.1244)&&9.09&&~(0.0396)&&68.10&&~(0.0297)&&19.99\\
\hline
\end{tabular}
\end{table}
\end{landscape}

\begin{landscape}
\begin{table}[!ht]
\centering
\scalefont{0.6}
\caption{Results of the analyzed methods for FSMD on the Golden instances.}
\label{GoldenFSMD}
\centering
\begin{tabular}{crccccccccccccccc}
\hline
  && \multicolumn{3}{c}{ILS-RVND \citep{Penna2013}} &&  \multicolumn{3}{c}{ILS-RVND-SP \citep{Subramanian2012}} &  \multicolumn{3}{c}{UHGS \citep{Vidal2014}} && \multicolumn{3}{c}{AILS}\\ 
\cline{3-5} \cline{7-9} \cline{11-13} \cline{15-17} 
Inst.&BKS&Avg~(\textit{gap})&Best&Time&&Avg~(\textit{gap})&Best&Time&&Avg~(\textit{gap})&Best&Time&&Avg~(\textit{gap})&Best&Time\\
\hline
\rowcolor{Gray}
3&$\mathbf{623.22^*}$&\textbf{623.22~(0.0000)}&\textbf{623.22}&4.58&&\textbf{623.22~(0.0000)}&\textbf{623.22}&\textbf{0.25}&&\textbf{623.22~(0.0000)}&\textbf{623.22}&10.20&&\textbf{623.22~(0.0000)}&\textbf{623.22}&8.58\\
4&$\mathbf{387.18^*}$&\textbf{387.18~(0.0000)}&\textbf{387.18}&2.85&&387.34~(0.0413)&\textbf{387.18}&0.23&&\textbf{387.18~(0.0000)}&\textbf{387.18}&11.40&&\textbf{387.18~(0.0000)}&\textbf{387.18}&\textbf{0.03}\\
\rowcolor{Gray}
5&$\mathbf{742.87^*}$&\textbf{742.87~(0.0000)}&\textbf{742.87}&5.53&&\textbf{742.87~(0.0000)}&\textbf{742.87}&0.22&&\textbf{742.87~(0.0000)}&\textbf{742.87}&12.00&&\textbf{742.87~(0.0000)}&\textbf{742.87}&\textbf{0.09}\\
6&$\mathbf{415.03^*}$&\textbf{415.03~(0.0000)}&\textbf{415.03}&3.37&&\textbf{415.03~(0.0000)}&\textbf{415.03}&\textbf{0.18}&&\textbf{415.03~(0.0000)}&\textbf{415.03}&13.20&&\textbf{415.03~(0.0000)}&\textbf{415.03}&9.45\\
\rowcolor{Gray}
13&$\mathbf{1491.86^*}$&1495.61~(0.2514)&\textbf{1491.86}&31.62&&1492.01~(0.0101)&\textbf{1491.86}&\textbf{1.91}&&\textbf{1491.86~(0.0000)}&\textbf{1491.86}&43.20&&1493.42~(0.1046)&\textbf{1491.86}&23.81\\
14&$\mathbf{603.21^*}$&\textbf{603.21~(0.0000)}&\textbf{603.21}&14.66&&605.00~(0.2967)&\textbf{603.21}&\textbf{1.61}&&\textbf{603.21~(0.0000)}&\textbf{603.21}&33.60&&\textbf{603.21~(0.0000)}&\textbf{603.21}&24.54\\
\rowcolor{Gray}
15&$\mathbf{999.82^*}$&1001.70~(0.1880)&\textbf{999.82}&15.33&&1001.03~(0.1210)&\textbf{999.82}&\textbf{1.47}&&\textbf{999.82~(0.0000)}&\textbf{999.82}&36.60&&1000.16~(0.0340)&\textbf{999.82}&31.58\\
16&$\mathbf{1131.00^*}$&1134.52~(0.3112)&\textbf{1131.00}&17.77&&1131.85~(0.0752)&\textbf{1131.00}&\textbf{1.44}&&\textbf{1131.00~(0.0000)}&\textbf{1131.00}&34.20&&1132.84~(0.1627)&\textbf{1131.00}&28.51\\
\rowcolor{Gray}
17&$\mathbf{1038.60^*}$&1041.12~(0.2426)&\textbf{1038.60}&49.18&&1042.48~(0.3736)&\textbf{1038.60}&6.39&&\textbf{1038.60~(0.0000)}&\textbf{1038.60}&68.40&&\textbf{1038.60~(0.0000)}&\textbf{1038.60}&\textbf{0.91}\\
18&$\mathbf{1800.80^*}$&1804.07~(0.1816)&\textbf{1800.80}&53.88&&1802.89~(0.1161)&\textbf{1800.80}&\textbf{4.75}&&1801.40~(0.0333)&1801.40&80.40&&\textbf{1800.80~(0.0000)}&\textbf{1800.80}&59.23\\
\rowcolor{Gray}
19&$1105.44$&1108.21~(0.2506)&\textbf{1105.44}&77.84&&1106.71~(0.1149)&\textbf{1105.44}&\textbf{10.62}&&1106.93~(0.1348)&\textbf{1105.44}&102.60&&\textbf{1105.44~(0.0000)}&\textbf{1105.44}&31.56\\
20&$\mathbf{1530.43^*}$&1540.32~(0.6462)&1530.52&88.02&&1534.23~(0.2483)&\textbf{1530.43}&\textbf{10.88}&&1531.82~(0.0908)&\textbf{1530.43}&168.00&&\textbf{1530.52~(0.0059)}&1530.52&50.35\\
\rowcolor{Gray}
Avg&&~(0.1726)&&30.39&&~(0.1164)&&3.33&&~(0.0216)&&51.15&&~(0.0256)&&22.39\\
\hline
\end{tabular}
\end{table}
\end{landscape}

\begin{table}[!ht]
\centering
\scalefont{0.55}
\caption{Results of the analyzed methods for FSMD on the instances introduced by \cite{Brandao2011}.}
\label{BrandaoFSMD}
\centering
\begin{threeparttable}
\begin{tabular}{ccccccccc}
\hline
&& \multicolumn{3}{c}{ILS-RVND-SP \citep{Subramanian2012}} && \multicolumn{3}{c}{AILS}\\ 
\cline{3-5} \cline{7-9} 
Inst.&BKS&Avg~(\textit{gap})&Best&Time&&Avg~(\textit{gap})&Best&Time\\
\hline
\rowcolor{Gray}
N1&$2211.63$\tnote{1}&\textbf{2219.66~(0.3631)}&\textbf{2212.77}&\textbf{39.60}&&2221.74~(0.4571)&2217.90&114.64\\
N2&$2811.37$\tnote{1}&2844.96~(1.1948)&2823.75&106.97&&\textbf{2821.60~(0.3639)}&\textbf{2812.28}&\textbf{89.50}\\
\rowcolor{Gray}
N3&$2234.57$&2234.85~(0.0125)&\textbf{2234.57}&\textbf{19.27}&&\textbf{2234.77~(0.0090)}&\textbf{2234.57}&59.41\\
N4&$1822.78$&1823.07~(0.0159)&\textbf{1822.78}&8.38&&\textbf{1822.78~(0.0000)}&\textbf{1822.78}&\textbf{7.25}\\
\rowcolor{Gray}
N5&$2016.79$&2019.26~(0.1225)&\textbf{2016.79}&\textbf{29.35}&&\textbf{2016.80~(0.0005)}&\textbf{2016.79}&73.62\\
Avg&&~(0.3418)&&40.71&&~(0.1661)&&68.88\\
\hline
\end{tabular}
\begin{tablenotes}
\item[1] Better BKS obtained by AILS.
\end{tablenotes}
\end{threeparttable}
\end{table}

\section{Final Remarks and Future Works}
\label{section:conclusao}

This paper has presented an Adaptive Iterated Local Search to solve the Heterogeneous Fleet Vehicle Routing Problem. AILS is an adaptive version of ILS with diversity control mechanisms tuned at runtime. This mechanism proved to be a very efficient hill-climbing strategy achieving robust results.

The introduced AILS solved five variants of the HVRP considering limited and unlimited fleets. This metaheuristic is simpler than the one presented by \cite{Maximo2021} as it does not have the hybridization with Path-Relinking (PR). In the computational experiments, we compared AILS with the state-of-the-art heuristic methods: ILS-RVND \citep{Penna2013}, ILS-RVND-SP \citep{Subramanian2012} and UHGS \citep{Vidal2014}. The experiments were carried out using the benchmark instances proposed by \cite{Golden1984} and \cite{Taillard1999} with up to 100 vertices. Considering this dataset, AILS had a better average result in the vast majority of the instances. The computational times required by AILS were reasonable when contrasting with the other algorithms. AILS was also compared with ILS-RVND-SP using the instances proposed by \cite{Li2007} and \cite{Brandao2011}. In these experiments, AILS found a better best-known solution on eight of the 15 instances. AILS had a better average solution than ILS-RVND-SP on all 15 instances. AILS showed a slower convergence and, for this reason, the computational times were lower only in four of the 15 instances.

As future work, we intend to evaluate the proposed algorithm for large HVRP instances and analyze the adjustments in AILS required for this setting. Another branch of study involves applying this approach to other HVRP variants, such as those considering multiple depots and time windows.

\section*{Acknowledgments}

The authors are grateful for the financial support provided by CNPq and FAPESP (2016/01860-1, 2019/22067-6). The Research was carried out using the computational resources of the Center for Mathematical Sciences Applied to Industry (CeMEAI) funded by FAPESP (grant 2013/07375-0).

\bibliographystyle{elsarticle-harv} \bibliography{Bibliography}

\end{document}